\input amstex
\documentstyle{amsppt}
\magnification=\magstep1
 \hsize 13cm \vsize 18.35cm \pageno=1
\loadbold \loadmsam
    \loadmsbm
    \UseAMSsymbols
\topmatter
\NoRunningHeads
\title A note on the $q$-analogue  of $p$-adic $\log$-gamma function
\endtitle
\author
  Taekyun Kim
\endauthor
 \keywords $p$-adic $q$-integrals, Euler numbers, q-Euler numbers,
 polynomials, sums of powers
\endkeywords

\abstract In this paper we prove that the $q$-analogue of Euler
numbers occur in the coefficients of some stirling type series for
the $p$-adic analytic $q$-$\log$-gamma function.

\endabstract
\endtopmatter

\document

{\bf\centerline {\S 1. Introduction}}

 \vskip 20pt

 Let $p$ be a fixed odd prime number. Throughout  this
paper $\Bbb Z$, $\Bbb Q$, $\Bbb Z_p$, $\Bbb Q_p$ and $\Bbb C_p$ will
respectively denote the ring of rational integers, the field of
rational numbers, the ring $p$-adic rational integers, the field of
$p$-adic rational numbers and the completion of the algebraic
closure of $\Bbb Q_p$. Let $v_p$ be the normalized exponential
valuation of $\Bbb C_p$ such that $|p|_p=p^{-v_p(p)}=p^{-1}$. If
$q\in\Bbb C_p$, we normally assume $|q-1|_p<p^{-\frac{1}{p-1}},$ so
that $q^x =\exp (x\log q)$ for $|x|_p\leq 1$. We use the notation
$$[x]_q=\frac{1-q^x}{1-q}, \text{ and  }
[x]_{-q}=\frac{1-(-q)^x}{1+q}.$$ Hence, $\lim_{q\rightarrow
1}[x]_q =1, $ for any $x$ with $|x|_p\leq 1$ in the present
$p$-adic case.

 For $d(= odd)$ a fixed positive
integer with $(p,d)=1$, let
$$\split
& X=X_d = \lim_{\overleftarrow{N} } \Bbb Z/ dp^N \Bbb Z , \text{  }
X_1 = \Bbb Z_p , \cr  & X^\ast = \underset {{0<a<d p}\atop
{(a,p)=1}}\to {\cup} (a+ dp \Bbb Z_p ), \cr & a+d p^N \Bbb Z_p =\{
x\in X | x \equiv a \pmod{dp^N}\},\endsplit$$ where $a\in \Bbb Z$
lies in $0\leq a < d p^N ,$ cf.[1-14].

In [3-7, 16], it is known that
$$\mu_{-q}(a+dp^N \Bbb Z_p)=(1+q)\frac{(-1)^aq^a}{1+q^{dp^N}}
=\frac{(-q)^a}{[dp^N]_{-q}},$$ is distribution on $X$ for $q\in \Bbb
C_p$ with $|1-q|_p< p^{-\frac{1}{p-1}}. $ This distribution yields
an integral as follows:
$$I_{-q}(f)=\int_{\Bbb Z_p} f(x)d\mu_{-q}(x)=\lim_{N\rightarrow
\infty}\frac{1}{[p^N]_{-q}}\sum_{x=0}^{p^N-1}f(x)(-q)^x, \text{ for
$f\in UD(\Bbb Z_p)$ }, \tag 1$$ which has a sense as we see readily
that the limit is convergent.

For $q=1$,  we have fermionic $p$-adic integral on $\Bbb Z_p$ as
follows:
$$I_{-1}=\int_{\Bbb
Z_p}f(x)d\mu_{-1}(x)=\lim_{N\rightarrow
\infty}\sum_{x=0}^{p^N-1}f(x)(-1)^x. $$ In view of notation,
$I_{-1}$ can be written symbolically as
$I_{-1}(f)=\lim_{q\rightarrow -1} I_{q}(f),$  where
$I_{q}(f)=\int_{\Bbb Z_p} f(x)d\mu_{q}(x)=\lim_{N\rightarrow
\infty}\frac{1}{[p^N]_{q}}\sum_{x=0}^{p^N-1}f(x)q^x ,$ see [3].

As the formula of the stirling asymptotic series, it was well known
that
$$ \log \left(\frac{\Gamma(x+1)}{\sqrt{2\pi}}\right)=(x-\frac{1}{2})\log x -x
+\sum_{n=1}^{\infty}\frac{(-1)^{n+1}}{n(n+1)}\frac{B_{n+1}}{x^n},
\text{ cf.[15]}, $$ where $B_n$ are called the $n$-th Bernoulli
numbers.

The purpose of this paper is to give  the new formula of the
$p$-adic $q$-analogue of
$\log\left(\frac{\Gamma(x+1)}{\sqrt{2\pi}}\right),$ which is related
to $q$-Euler numbers. That is, we prove that the $q$-analogue of
Euler numbers occur in the coefficients of some stirling type series
for $p$-adic analytic $q$-$\log$-gamma functions.
 \vskip 20pt

{\bf\centerline {\S 2. $p$-adic $q$-$\log$-gamma function}}

 \vskip 20pt

Let us include some remarks about the factorial function, we define
$0!=1$ and may compute further values by the relation
$(n+1)!=n!(n+1)$. For large $n$ the function is very large. A
convenient approximation for large $n$ is the stirling formula:
 $$n!\sim\sqrt{2\pi n}(\frac {n}{e})^n,
\text{ $(e=2.718\cdots)$, cf.[15]}, \tag1-2$$ where $\sim$ means
that the ratio of two sides of (1-2) approaches 1 as $n$ approaches
infinity.

From (1-2) we can derive $$ \log
\left(\Gamma(x+1)/\sqrt{2\pi}\right)=(x+B_1)\log x -x
+\sum_{n=1}^{\infty}\frac{(-1)^{n+1}}{n(n+1)}\frac{B_{n+1}}{x^n},
\text{ cf. [4, 15]}, \tag2$$ where $B_n$ are called the $n$-th
Bernoulli numbers.

For any non-negative integer $m$, we define the $q$-Euler
polynomials as follows:
$$\int_{\Bbb Z_p}[x+y]_q^m
d\mu_{-q}(y)=E_{m,q}(x)=[2]_q\left(\frac{1}{1-q}\right)^m\sum_{i=0}^m\binom{m}{i}(-1)^i\frac{q^x}{1+q^{i+1}}.\tag3$$
From (3), we can also derive $q$-Euler numbers, $E_{n,q}$,  as
$E_{n,q}(0)=E_{n,q}.$ Note that $\lim_{q\rightarrow 1}E_{n,q}=E_n$,
where $E_n$ are ordinary Euler numbers which are defined by
$\frac{2}{e^t+1}=\sum_{n=0}^{\infty}E_n\frac{t^n}{n!}. $ By the
simple calculation, it is easy to show that
$$\left((1+x)\log(1+x)\right)^{\prime}=1+\log({1+x}) =1+\sum_{n=1}^{\infty}\frac{(-1)^{n+1}}{n(n+1)}x^n, \tag4 $$
where
$\left((1+x)\log(1+x)\right)^{\prime}=\frac{d}{dx}\left((1+x)\log(1+x)\right).$

From (4) we derive
$$(1+x)\log(1+x)=\sum_{n=1}^{\infty}\frac{(-1)^{n+1}}{n(n+1)}x^{n+1}+x+c, \text{ where $c$ is constant}.\tag5$$
If we take $x=0$, then we have $c=0$. By (3) and (4), we easily see
that
$$(1+x)\log(1+x)=\sum_{n=1}^{\infty}\frac{(-1)^{n+1}}{n(n+1)}x^{n+1}+x.\tag6$$
We now consider $p$-adic locally analytic function $G_{p,q}(x)$ on
$\Bbb C_p \backslash\Bbb Z_p$ by

$$G_{p,q}(x)=\int_{\Bbb Z_p}[x+z]_q\left(
\log[x+z]_q -1\right)d\mu_{-q}(z). \tag 7$$ From (1) we can easily
derive
$$qI_{-q}(f_1)+I_{-q}(f)=[2]_qf(0), \text{ where $f_1$ is
translation with $f_{1}(x)=f(x+1)$.}\tag8$$ By (7) and (8), we
easily see that
$$qG_{p,q}(x+1)+G_{p,q}(x)=[2]_q\left([x]_q(\log[x]_q-1)\right).$$

It is easy to see that
$$[x+z]_q=\frac{1-q^{x+z}}{1-q}=\frac{1-q^x
+q^x(1-q^z)}{1-q}=[x]_q+q^x[z]_q. \tag9$$ By(6) and (9) we see that
$$\split &[x+z]_q(\log[x+z]_q-1)\cr
&=[z]_q+[x]_q\sum_{n=1}^{\infty}\frac{(-q^x)^{n+1}}{n(n+1)}\frac{[z]_q^{n+1}}{[x]_q^{n+1}}
+([x]_q+q^x[z]_q)\log[x]_q-([x]_q+[z]_q).\endsplit\tag10$$ From (3),
(7) and (10), we note that
$$G_{p,q}(x)=([x]_q+q^xE_{1,q})\log[x]_q-[x]_q
+\sum_{n=1}^{\infty}\frac{(-q^x)^{n+1}}{n(n+1)}\frac{1}{[x]_q^n}E_{n+1,q}.$$

Therefore we obtain the following:

\proclaim{ Theorem A} For $x\in\Bbb C_p\backslash\Bbb Z_p $, we have
$$G_{p,q}(x)=([x]_q-q^x\frac{1}{[2]_{q^2}})\log[x]_q-[x]_q
+\sum_{n=1}^{\infty}\frac{(-q^x)^{n+1}}{n(n+1)}\frac{1}{[x]_q^n}E_{n+1,q},\tag
11$$ and
$$qG_{p,q}(x+1)+G_{p,q}(x)=[2]_q\left([x]_q(\log[x]_q-1)\right).\tag12$$
\endproclaim

Remark. The above Theorem A seems to be the $p$-adic $q$-analogue of
$\log \frac{\Gamma(x+1)}{\sqrt{2\pi}},$ which is related to
$q$-Euler numbers. In [4], $q$-Bernoulli numbers defined by
$$\int_{\Bbb Z_p} q^{-x}[x]_q^nd\mu_q(x)=\beta_{n,q}.$$
For $x\in\Bbb C_p \backslash \Bbb Z_p$, we consider the $p$-adic
$q$-$\log$-gamma function as follows:
$$T_{p,q}(x)=\int_{\Bbb Z_p}
q^{-y-x}[x+y]_q\left(\log[x+y]_q-1\right)d\mu_q(y). \tag13$$ From
(13) and (6) it seems to be derived the following interesting
formula:
$$
T_{p,q}(x)=(q^{-x}[x]_q\beta_{0,q}+\beta_{1,q})\log[x]_q-q^{-x}[x]_q\beta_{0,q}
+\sum_{n=1}^{\infty}\frac{(-1)^{n+1}q^{nx}}{n(n+1)}\frac{\beta_{n+1,q}}{[x]_q^n}.$$

\vskip 20pt

\quad Taekyun Kim

\quad EECS, Kyungpook National University, Taegu 702-701, S. Korea

\quad e-mail:\text{ tkim$\@$knu.ac.kr; tkim64$\@$hanmail.net}


\Refs \widestnumber\key{999999}
\ref \key 1
 \by M. Cenkci, M. Can, V. Kurt
  \paper $p$-adic interpolation functions and Kummer-type congruences for $q$-twisted and
  $q$-generalized twisted Euler numbers
 \jour   Adv. Stud. Contemp. Math.
 \yr 2004
\pages 203-216 \vol 9 \endref


\ref \key 2
 \by M. Cenkci, M. Can
  \paper Some results on $q$-analogue of the Lerch zeta function
 \jour   Adv. Stud. Contemp. Math.
 \yr 2006
\pages 213-223 \vol 12 \endref



\ref \key 3
 \by  T. Kim
  \paper  On p-adic interpolating function for q-Euler numbers and its derivatives
 \jour J. Math. Anal. Appl.
 \yr 2008
\pages 598-608 \vol 339 \endref


\ref \key 4
 \by  T. Kim
  \paper  $q-$Volkenborn integration
 \jour  Russ. J. Math. Phys.
 \yr 2002
\pages 288--299 \vol 9 \endref


\ref \key 5
 \by  T. Kim
  \paper  A note on $p$-adic invariant integral in the rings of
  $p$-adic integers
 \jour  Advan. Stud. Contemp. Math.
 \yr 2006
\pages 95--99 \vol 13 \endref

\ref \key 6
 \by  T. Kim
  \paper  q-Extension of the Euler formula and trigonometric functions
 \jour  Russ. J. Math. Phys.
 \yr 2007
\pages 275-278 \vol 14 \endref

\ref \key 7
 \by  T. Kim
  \paper  A new approach to q-zeta functions
 \jour Journal of Computational Analysis and Applications
 \yr 2007
\pages 395-400 \vol 9 \endref

\ref \key 8
 \by  L.C. Jang, S.D. Kim, H. K. Park, Y.S. Ro
  \paper  A note on Euler number and polynomials
 \jour J. Inequal. Appl.
 \yr 2006
\pages Art. ID 34602, 5 pp \vol 2006 \endref

\ref \key 9
 \by  H. Ozden, Y. Simsek, S.-H. Rim, I. N.  Cangul
  \paper  A note on $p$-adic $q$-Euler measure
 \jour Adv. Stud. Contemp. Math.
 \yr 2007
\pages 233-239 \vol 14 \endref

\ref \key 10
 \by  C. S. Ryoo
  \paper  A note on q-Bernoulli numbers and polynomials
 \jour Applied Mathematics Letters
 \yr 2007
\pages 524-531 \vol 20 \endref

\ref \key 11
 \by  C. S. Ryoo
  \paper  A numerical computation on the structure of the roots of q-extension of Genocchi
  polynomials
 \jour Applied Mathematics Letters
 \yr 2007
\pages doi:10.1016/j.aml.2007.05.005   \vol  \endref

\ref \key 12
 \by  Y. Simsek
  \paper  Twisted $(h,q)$-Bernoulli numbers and polynomials related to twisted $(h,q)$-zeta function and
  $L$-function
 \jour J. Math. Anal. Appl.
 \yr 2006
\pages 790--804  \vol 324 \endref

\ref \key 13
 \by  Y. Simsek
  \paper  On $p$-adic twisted $q\text{-}L$-functions related to generalized twisted Bernoulli numbers
 \jour Russ. J. Math. Phys.
 \yr 2006
\pages 340--348  \vol 13 \endref

\ref \key 14
 \by  Y. Simsek
  \paper On twisted $q$-Hurwitz zeta function and $q$-two-variable $L$-function
 \jour  Appl. Math. Comput.
 \yr 2007
\pages 466--473  \vol 187 \endref

\ref \key 15
 \by  D. Zill, M. R. Cullen
  \book Advanced Engineering Mathematics
 \publ  Jones and Bartlett
 \yr 2005
\pages   \endref

\ref \key 16
 \by  T. Kim, J. Y. Choi, J. Y. Sug,
  \paper Extended $q$-Euler numbers and polynomials associated with fermionic $p$-adic $q$-integral on $\Bbb Z\sb p$
 \jour  Russ. J. Math. Phys.
 \yr 2007
\pages 160--163  \vol 14 \endref



\endRefs
\enddocument